\documentclass[graybox]{svmult}

\usepackage{helvet}         
\usepackage{courier}        
\usepackage{type1cm}        
\usepackage{makeidx}         
\usepackage{graphicx}        
\usepackage{multicol}        
\usepackage[bottom]{footmisc}

\usepackage{amsmath}
\usepackage{amsxtra}
\usepackage{amscd}
\usepackage{amsfonts}
\usepackage{amssymb}

\newcommand{\Glie}{\mathfrak{g}}             
\newcommand{\BC}{\mathbb{C}}            
\newcommand{\BZ}{\mathbb{Z}}            

\newcommand{\BV}{\mathbf{V}}             
\newcommand{\BQ}{\mathbf{Q}}             
\newcommand{\BP}{\mathbf{P}}             
\newcommand{\Hom}{\mathrm{End}}          
\newcommand{\Id}{\mathrm{Id}}            
\newcommand{\Sm}{\mathbb{S}}             

\newcommand{\super}{\mathbb{Z}_2}        
\newcommand{\even}{\overline{0}}         
\newcommand{\odd}{\overline{1}}          

\begin{document}

\title*{Two-parameter quantum general linear supergroups}
\author{Huafeng Zhang}
\institute{Laboratoire Paul Painlev\'e \& 
Universit\'e de Lille, 
59655 Villeneuve d'Ascq, France, \email{Huafeng.Zhang@math.univ-lille1.fr}
}
\maketitle

\abstract{The universal R-matrix of two-parameter quantum general linear supergroups is computed explicitly based on the RTT realization of Faddeev--Reshetikhin--Takhtajan.}

\section*{Introduction}
Fix $r,s$ non-zero complex numbers whose ratio $\frac{r}{s}$ is not a root of unity. Let $M,N$ be positive integers and $\mathfrak{g} := \mathfrak{gl}(M,N)$ be the general linear Lie superalgebra. The  enveloping algebra  $U(\mathfrak{g})$ as a Hopf superalgebra admits a two-parameter deformation $U_{r,s}(\mathfrak{g})$ which is neither commutative nor cocommutative. In this paper we compute its {\it universal R-matrix}, an invertible element in a completed tensor square $\mathcal{R} \in U_{r,s}(\mathfrak{g})^{\widehat{\otimes} 2}$ satisfying 
\begin{displaymath}
\Delta^{\mathrm{cop}}(x) = \mathcal{R} \Delta(x) \mathcal{R}^{-1} \quad \mathrm{for}\ x \in U_{r,s}(\mathfrak{g}),
\end{displaymath} 
together with other favorable properties.
In the non-graded case $N = 0$, Benkart--Witherspoon \cite{BW1,BW2} proved the existence of universal R-matrix, and derived from it  a braided structure in the category of finite-dimensional representations; the exact formula of universal R-matrix was unknown. Recently it was shown \cite{JL} that $U_{r,s}(\mathfrak{gl}(M))$ can be recovered from a special R-matrix in the spirit of Faddeev--Reshetikhin--Takhtajan \cite{FRT}, the RTT realization.  

In this paper we define the two-parameter quantum supergroup $U_{r,s}(\Glie)$ by RTT realization, based on a suitable R-matrix on the vector superspace $\BC^{M|N}$. Our main result, Equations \eqref{equ: universal R matrix factorization}--\eqref{equ: partial R matrix}, is a factorization formula for the universal R-matrix $\mathcal{R}$ in terms of RTT generators. (This idea was previously applied to
 the quantum affine superalgebra of $\mathfrak{gl}(1,1)$; see \cite{Z2}.)

Let us compare with earlier works on universal R-matrices: \cite{Rosso} for $U_q(\mathfrak{sl}_M)$; \cite{Y1} for $U_q(\mathfrak{gl}(M,N))$; \cite{AY} for (quantum doubles of) Nichols algebras, which are believed to include two-parameter quantum (super)groups. In these works a key step is to construct root vectors by Lusztig isomorphisms or q-brackets. In our approach the root vectors are already encoded in the definition of the algebra. There is another two-parameter quantum supergroup $U_{q_1,q_2}(\mathfrak{sl}(2,1))$ proposed by R.B. Zhang \cite{RB}: for $q_1= q_2 =q$ it is $U_q(\mathfrak{sl}(2,1))$, while for $q_1 \neq q_2$ its comultiplication is not yet clear.

\section{RTT realization and orthogonality}  \label{sec: YBE}
We define $U_{r,s}(\mathfrak{gl}(M,N))$ following Faddeev--Reshetikhin--Takhtajan \cite{FRT}, and prove an orthogonality property for the associated Hopf pairing.

Let $\BV = \BC^{M|N}$ be the vector superspace with basis $(v_i)_{1\leq i \leq M+N}$ and parity:
$ |v_i| = |i| = \even$ if $i \leq M$ and $|v_i| = |i| = \odd$ if $i > M$. Define the elementary matrices $E_{ij} \in \Hom \BV: \ v_k \mapsto \delta_{jk} v_i$.
 Define the two-parameter Perk--Schultz matrix $R \in \Hom(\BV^{\otimes 2})$ by
\begin{equation}   \label{equ: R formula}
(r\sum_{i \leq M} + s\sum_{i > M}) E_{ii} \otimes E_{ii} + (\sum_{i>j} + rs \sum_{i<j}) E_{ii} \otimes E_{jj} + (r-s) \sum_{i<j} (-1)^{|i|} E_{ji} \otimes E_{ij}.
\end{equation}
Recall the super tensor product. For $V = V_{\even} \oplus V_{\odd}$ a  vector superspace  and $p \in \super = \{\even,\odd\}$, let $(\Hom V)_p$ denote the set of linear endomorphisms $g\in \Hom V$ such that $g(V_q) \subseteq V_{p+q}$ for all $q \in \super$. This makes $\Hom V$ a superalgebra. Let $W$ be another vector superspace.  For $f \in \Hom W$ and $g \in (\Hom V)_p$  the super tensor product $f \otimes g \in \Hom(W\otimes V)$ is defined by 
$$ f \otimes g: w \otimes v \mapsto (-1)^{p q} f(w) \otimes g(v) \quad \mathrm{for}\ w \in W_q\ \mathrm{and}\ v \in V. $$
If $V,W$ are finite-dimensional, this identifies the tensor product superalgebra $\Hom W \otimes \Hom V$ with $\Hom (W\otimes V)$. Let us define three elements of $\Hom (\BV^{\otimes 3})$:
$$ R_{12} = R \otimes 1,\quad R_{23} = 1 \otimes R, \quad R_{13} = (c_{\BV,\BV} \otimes 1) R_{23} (c_{\BV,\BV} \otimes 1). $$
Here $c_{\BV,\BV} \in \Hom (\BV^{\otimes 2}): v_i \otimes v_j \mapsto (-1)^{|i||j|} v_j \otimes v_i$ is the graded flip.

\begin{lemma}[Yang--Baxter Equation] \label{lem: QYBE}
$R_{12}R_{13}R_{23} = R_{23}R_{13}R_{12}$. 
\end{lemma}
\begin{proof}
Set $\widehat{R} := c_{\BV,\BV} R \in \Hom(\BV^{\otimes 2})$. Define $\widehat{R}_{12}$ and $\widehat{R}_{23} \in \Hom (\BV^{\otimes 3})$ in the obvious way. The Yang--Baxter equation is equivalent to the braid relation \begin{equation} \label{braid}
\widehat{R}_{12}\widehat{R}_{23}\widehat{R}_{12} = \widehat{R}_{23} \widehat{R}_{12} \widehat{R}_{23} \in \mathrm{End}(\BV^{\otimes 3}).
\end{equation}
To indicate the dependence on $r,s,M,N$, we shall also let $\widehat{R}(r,s,M,N)$ denote $\widehat{R}$. Since $\widehat{R}$ is of even parity, the validity of Equation \eqref{braid} is independent of the $\super$-grading on $\BV = \BC^{M|N}$. 
 Observe that
$ \widehat{R}(r,s,M,0) = \widehat{R}(s,r,0,M)$.
By \cite[Proposition 5.5]{BW2}, Equation \eqref{braid} holds for the matrix $\widehat{R}(r,s,n,0)$. So it holds for $\widehat{R}(r,s,M+N,0) =: \mathcal{S}$ and $\widehat{R}(r,s,0,M+N) =: \mathcal{S}'$. After ignoring the super structure, the vector spaces $\BC^{M+N|0}, \BC^{0|M+N}$ and $\BC^{M|N}$ are the same. So we view $\mathcal{S}, \mathcal{S}' \in \Hom(\BV^{\otimes 2})$. 

We prove that \eqref{braid} applied to $v_a \otimes v_b \otimes v_c$ is true if $a,b,c \in \{1,2,\cdots,M+N\}$ are two-by-two distinct. Let $\mathcal{S}_{kl}^{ij}$ be the coefficient of $v_k \otimes v_l$ in the vector $\mathcal{S}(v_i \otimes v_j)$. Then $\mathcal{S}_{kl}^{ij} \neq 0$ implies $\{i,j\} = \{k,l\}$, and for $i \neq j$ we have 
\begin{equation} \label{RS}
\widehat{R}(v_i \otimes v_j) = \mathcal{S}_{ij}^{ij} v_i \otimes v_j + (-1)^{|i||j|} \mathcal{S}_{ji}^{ij} v_j \otimes v_i. 
\end{equation}
Apply $\mathcal{S}_{12} \mathcal{S}_{23} \mathcal{S}_{12} = \mathcal{S}_{23}\mathcal{S}_{12}\mathcal{S}_{23}$ to $v_a \otimes v_b \otimes v_c$, and let $C_{ijk}$ be the coefficient of $v_i \otimes v_j \otimes v_k$. Then $C_{ijk} \neq 0$ only if $ijk$ is a permutation of $abc$. Based on the relation \eqref{RS} of $\widehat{R}$ and $\mathcal{S}$, one proves that
$$\widehat{R}_{12}\widehat{R}_{23}\widehat{R}_{12}(v_a\otimes v_b \otimes v_c) = \sum_{ijk} s_{ijk} C_{ijk} v_i \otimes v_j \otimes v_k = \widehat{R}_{23} \widehat{R}_{12} \widehat{R}_{23}(v_a \otimes v_b \otimes v_c)  $$
where $s_{ijk} = \pm$ is a signature depending on the permutation $ijk$ of $abc$. 

Based on the braid relations on $\mathcal{S}$ and $\mathcal{S}'$, one shows that $\eqref{braid}$ applied to $v_a \otimes v_b \otimes v_c$ is true if $abc$ is a permutation of $iij$ such that $i \leq M $ or $i,j > M$. 

We are reduced to the case $M=N=1$ and to show that the braid relation applied to $v_1 \otimes v_2 \otimes v_2,\ v_2 \otimes v_1 \otimes v_2,\ v_2 \otimes v_2 \otimes v_1$ holds. Set $\mathcal{T} := \widehat{R}(r,s,1,1)$. Consider the second vector $u := v_2 \otimes v_1 \otimes v_2$ as an example:
\begin{align*}
\mathcal{T}_{12}\mathcal{T}_{23}T_{12}(u) &= \mathcal{T}_{12}\mathcal{T}_{23}(v_1 \otimes v_2 \otimes v_2) = -s \mathcal{T}_{12}(v_1 \otimes v_2 \otimes v_2) \\
&= -s(r-s) v_1 \otimes v_2 \otimes v_2 - rs^2 v_2 \otimes v_1 \otimes v_2 \\
&= \mathcal{T}_{23}((r-s)v_1 \otimes v_2 \otimes v_2 - rs^2 v_2 \otimes v_2 \otimes v_1)  \\
&= \mathcal{T}_{23}\mathcal{T}_{12}((r-s)v_2 \otimes v_1 \otimes v_2 + rs v_2 \otimes v_2 \otimes v_1) = \mathcal{T}_{23}\mathcal{T}_{12}\mathcal{T}_{23}(u).
\end{align*}
The first and the third vectors can be checked in the same way.  \hfill $\Box$
\end{proof}

\begin{definition}   \label{definition: two-parameter quantum supergroup}
$U := U_{r,s}(\mathfrak{gl}(M,N))$ is the superalgebra generated by the coefficients of matrices $T = \sum_{i\leq j} t_{ji} \otimes E_{ji},\ S = \sum_{i\leq j}s_{ij} \otimes E_{ij} \in U \otimes \Hom \BV$ of even parity (so that $s_{ij}$ and $t_{ji}$ are of parity $|i|+|j|$) with relations
\begin{eqnarray*}
 R_{23}T_{12}T_{13} = T_{13}T_{12}R_{23},\ R_{23}S_{12}S_{13} = S_{13}S_{12}R_{23},\  R_{23}T_{12}S_{13} = S_{13}T_{12}R_{23},
\end{eqnarray*}
and the $s_{ii},t_{ii}$ are invertible for $1 \leq i \leq M+N$.  
\end{definition}
$U$ is a Hopf superalgebra with coproduct $\Delta$ and counit $\varepsilon$:
$$ \Delta(s_{ij}) = \sum_{k} s_{ik} \otimes s_{kj},\quad \Delta(t_{ji}) = \sum_k t_{jk} \otimes t_{ki},\quad \varepsilon(s_{ij}) = \varepsilon(t_{ji}) = \delta_{ij}. $$
The antipode $\Sm: U \longrightarrow U$ is an anti-automorphism of superalgebra defined by equations $(\Sm \otimes \mathrm{Id})(S) = S^{-1}, (\Sm \otimes \mathrm{Id})(T) = T^{-1}$ in $U \otimes \Hom \BV$. Let $U^+$ (resp. $U^-$) be the subalgebra of $U$ generated by the $s_{ij}, s_{kk}^{-1}$ (resp. the $t_{ji}, t_{kk}^{-1}$) for $i \leq j$; these are sub-Hopf-superalgebras. 
Algebra $U$ is graded by the weight lattice $\BP := \oplus_{i=1}^{M+N} \BZ \epsilon_i$; we set $s_{ij}$ and $t_{ji}$ to be of weight $\pm (\epsilon_i-\epsilon_j)$ respectively. The weight grading restricts to subalgebras $U^{\pm}$.

We interpret  Definition \ref{definition: two-parameter quantum supergroup} as a quantum double construction, following \cite[\S 3.1.3]{Z1}. There exists a unique bilinear form $\varphi: U^+ \times U^- \longrightarrow \BC$ such that
\begin{equation} \label{equ: R Hopf pairing}
\sum_{ijkl} \varphi(s_{ij},t_{kl}) E_{kl} \otimes E_{ij} = R \in (\Hom \BV)^{\otimes 2},
\end{equation}
and for $a,a' \in U^+$ and $b,b' \in U^-$ super homogeneous
\begin{gather*}
\varphi(a,b b') = \varphi_2(\Delta(a),b\otimes b'),\quad \varphi(aa',b) = (-1)^{|a||a'|} \varphi_2(a' \otimes a, \Delta(b)).
\end{gather*}
Here $\varphi_2(a \otimes a', b\otimes b') = (-1)^{|a'||b|} \varphi(a,b)\varphi(a',b')$.
 Such a form is called {\it Hopf pairing}. The quantum double $U^+ \otimes U^-$ is isomorphic to $U$ as Hopf superalgebras via the multiplication map. This implies that in $U$:
\begin{equation}  \label{rem: Hopf pairing}
b a  =  (-1)^{|a_{(1)}||b| + (|b_{(2)}|+|b_{(3)}|)|a_{(2)}| + |a_{(3)}||b_{(3)}|} \varphi (a_{(1)}, \Sm (b_{(1)})) a_{(2)} b_{(2)} \varphi(a_{(3)}, b_{(3)}).
\end{equation}
Here $a_{(1)} \otimes a_{(2)} \otimes a_{(3)} = (\Delta \otimes \Id)\Delta(a)$ is the Sweedler notation.

The Hopf pairing respects the weight grading: for $x \in U^+$ and $y \in U^-$ being of weight $\alpha$ and $\beta$ respectively, $\varphi(x,y) \neq 0$ only if $\alpha+\beta = 0$.
 
Let $\tau: \Hom \BV \longrightarrow \Hom \BV$ be the transposition $E_{ij} \longrightarrow (-1)^{|i|+|i||j|} E_{ji}$. Lemma \eqref{lem: QYBE} affords  a vector representation $\rho$ of $U$ on $\BV$:
\begin{equation} \label{equ: vector repre}
(\rho \otimes 1)(S) = (\tau \otimes 1)(R),\quad (\rho \otimes 1)(T) = rs (\tau \otimes 1)(c_{\BV,\BV}R^{-1}c_{\BV,\BV}).
\end{equation}

\begin{lemma} \label{cor: weight}
Let $1\leq i,j,k \leq M+N$ be such that $j \leq k$. Then
\begin{align*}
s_{ii}s_{jk} &= \varphi(s_{ii},t_{jj})\varphi(s_{ii},t_{kk})^{-1} s_{jk}s_{ii},\quad t_{ii}s_{jk} = \varphi(s_{jj},t_{ii})^{-1}\varphi(s_{kk},t_{ii})s_{jk}t_{ii}, \\
t_{ii}t_{kj} &= \varphi(s_{jj},t_{ii})\varphi(s_{kk},t_{ii})^{-1} t_{kj}t_{ii},\quad s_{ii} t_{kj} =  \varphi(s_{ii},t_{jj})^{-1}\varphi(s_{ii},t_{kk}) t_{kj}s_{ii}.
\end{align*}
\end{lemma}
\proof
For the second identity, by Equation \eqref{rem: Hopf pairing}
\begin{displaymath}
t_{ii} s_{jk} = \varphi(s_{jj},\Sm(t_{ii})) s_{jk}t_{ii} \varphi(s_{kk},t_{ii}) = \varphi(s_{jj},t_{ii})^{-1}\varphi(s_{kk},t_{ii})s_{jk}t_{ii}.
\end{displaymath}
Here we have used the three-fold coproduct formula of $t_{ii},s_{jk}$, and the fact that $\varphi(s_{ab},t_{ii}) = 0$ if $a < b$. The fourth identity can be proved similarly. 

For the first identity, by comparing the coefficients of $v_i \otimes v_j$ in the identical vectors $R_{23}S_{12}S_{13}(v_i\otimes v_k) = S_{13}S_{12} R_{23}(v_i\otimes v_k) \in U \otimes \BV^{\otimes 2}$ we obtain
\begin{displaymath}
x s_{ii}s_{jk} + y s_{ji}s_{ik} = z s_{jk}s_{ii} + w s_{ji}s_{ik}
\end{displaymath}
for certain $x,z \in \{1,rs,r,s\}$ and $y,w \in \{0,r-s,s-r\}$. Here we set $s_{pq} = 0$ if $p > q$. We prove that $s_{ii}s_{jk} \in \BC s_{jk}s_{ii}$. If not, then $j<i<k$, in which case $y = (-1)^{|i|}(r-s) = w$ and $x s_{ii} s_{jk} = z s_{jk}s_{ii}$, a contradiction. Now the first identity is obtained from the vector representation \eqref{equ: vector repre}:
$$\rho(s_{ii}) = \sum_k \varphi(s_{ii},t_{kk}) E_{kk}, \quad \rho(s_{jk}) = \varphi(s_{jk},t_{kj}) (-1)^{|k|+|j|} E_{jk} \quad \textrm{for}\ j < k.$$
The third identity can be proved in the same way. 
\hfill $\Box$

It follows that  a vector $x \in U$ is of weight $\sum_i \lambda_i \epsilon_i$ if  and only if 
\begin{equation*}  
s_{ii} x s_{ii}^{-1} = \varphi(s_{ii},t_{ii})^{\lambda_i} (rs)^{\sum_{j<i}\lambda_j} x,\quad t_{ii}xt_{ii}^{-1} = \varphi(s_{ii},t_{ii})^{-\lambda_i} (rs)^{\sum_{j\leq i}\lambda_j} x.
\end{equation*}

Let us define the modified RTT generators
\begin{equation}  \label{equ: modified RTT}
a_{ij} := s_{ii}^{-1}s_{ij},\quad b_{ji} := t_{ji}t_{ii}^{-1} \quad \mathrm{for}\ 1 \leq i < j \leq M+N. 
\end{equation}
The $a_{ij}$ form a subset $X$ and generate a subalgebra $U^>$ of $U^+$. Similarly, the $b_{ji}$ form a subset $Y$ and generate a subalgebra $U^<$ of $U^-$. Let $H^+$ (resp. $H^-$) be the subalgebra of $U^>$ (resp. $U^<$) generated by the $s_{ii}$ (resp. the $t_{ii}$).

 $X, Y$ are totally ordered sets with lexicographic ordering:
$a_{ij}\prec a_{kl}$ and $b_{ji} \prec b_{lk}$ if either $(i < k)$ or $(i = k,j < l)$.

\begin{lemma}   \label{lem: technical lemma}
Fix $1\leq i < j \leq M+N$ and $p \in \BZ_{>0}$. Let $x_1,x_2,\cdots,x_p \in X$ and $y_1,y_2,\cdots,y_p \in Y$ be such that $x_l \succeq a_{ij}$ and $y_l \succeq b_{ji}$ for all $1 \leq l \leq p$. 
\begin{enumerate}
\item[(A)] We have $\varphi(a_{ij},b_{ji}) = (-1)^{|i|}(s^{-1}-r^{-1})$. 
\item[(B)] If $\varphi(a_{ij},y_1y_2\cdots y_p) \neq 0$, then $p = 1, y_1 = b_{ji}$.
\item[(C)] If $\varphi(x_1x_2\cdots x_p, b_{ji}) \neq 0$, then $p = 1$ and $x_1 = a_{ij}$.
\end{enumerate}
\end{lemma}
\proof
Let us first prove an auxiliary result:
\begin{enumerate}
\item[(D)] If $a \in U^>, b \in U^<$ and $x_{\pm} \in H^{\pm}$, then $\varphi(x_+a, x_-b) = \varphi(x_+,x_-)\varphi(a,b)$. 
\end{enumerate} 

One may assume that $x_+$ is a product of the $s_{ii}^{\pm 1}$ so that $\varphi(x_+,1) = 1$ and $\Delta(x^+) = x^+ \otimes x^+$. By definition, $\Delta(x_+a) - x_+ \otimes x_+a$ is a sum of $x_i' \otimes y_i'$ where each $x_i'$ is of non-zero weight and so $\varphi(x_i',x_-) = 0$. By Equation \eqref{rem: Hopf pairing},
\begin{displaymath}
\varphi(x_+a, x_-b) = \varphi_2(x_+\otimes x_+a, x_- \otimes b) = \varphi(x_+,x_-)\varphi(x_+a,b).
\end{displaymath}
$\Delta(b) - b \otimes 1$ is a sum of $x_i'' \otimes y_i''$ where each $y_i''$ is of non-zero weight and so $\varphi(x_+,y_i'') = 0$. This implies
\begin{displaymath}
\varphi(x_+a,b) =  \varphi_2(a \otimes x_+, b\otimes 1) = \varphi(a,b)\varphi(x_+,1) = \varphi(a,b).
\end{displaymath}
 This proves (D). We are able to compute $\varphi(a_{ij},b_{ji})$:
\begin{align*}
 \varphi(s_{ij},t_{ji}) &= \varphi(s_{ii}a_{ij},b_{ji}t_{ii}) = \varphi(s_{ii}a_{ij},\varphi(s_{ii},t_{ii})^{-1}\varphi(s_{jj},t_{ii})t_{ii}b_{ji}) \\
 &= \varphi(s_{jj},t_{ii})\varphi(a_{ij},b_{ji}) = rs \varphi(a_{ij},b_{ji}).
\end{align*}
(A) follows from Equations \eqref{equ: R formula} and \eqref{equ: R Hopf pairing}.
For (B), the first tensor factors in $\Delta(a_{ij}) - a_{ij} \otimes s_{ii}^{-1}s_{jj}$, being either 1 or $x \in X$ with $x \prec a_{ij}$, are orthogonal to $y_1\succeq b_{ji}$. So $\varphi(a_{ij},y_1y_2\cdots y_p) = \varphi(a_{ij},y_1) \varphi(s_{ii}^{-1}s_{jj},y_2\cdots y_p)$. Now $\varphi(a_{ij},y_1) \neq 0$ forces $p = 1$ and $y_1 = b_{ji}$. (C) is proved similarly.  \hfill $\Box$

\begin{lemma}  \label{cor: technical corollary}
Fix $1\leq i < j \leq M+N$. Let $x_1,x_2,\cdots,x_p \in \{x\in X\ |\ x \succ a_{ij}\}$ and $y_1,y_2,\cdots,y_q \in \{y \in Y\ |\ y \succ b_{ji}\}$. Let $m,n \in \BZ_{\geq 0}$. Then 
\begin{align*}
\varphi(x_1x_2\cdots x_p a_{ij}^m, y_1y_2\cdots y_q b_{ji}^n) &= \varphi_2(x_1x_2\cdots x_p \otimes a_{ij}^m, y_1y_2\cdots y_q \otimes b_{ji}^n), \\
\varphi(a_{ij}^m,b_{ji}^n) &=  \delta_{mn} (m)_{\tau_{ij}}^! \varphi(a_{ij},b_{ji})^m. 
\end{align*}
Here $(m)_{u} := \prod_{k=1}^m \frac{u^{k}-1}{u-1}$ and $\tau_{ij} := (-1)^{|i|+|j|} (rs)^{-1}\varphi(s_{ii},t_{ii})\varphi(s_{jj},t_{jj})$.
\end{lemma}
\proof
By induction on $\max(m,n)$: the case $m = n = 0$ is trivial. Assume $m > 0$ (the case $n>0$ can be treated similarly). The left hand side of the first formula becomes (we set $\theta_1 := |a_{ij}||a_{ij}^{m-1}x_1x_2\cdots x_p|$)
\begin{displaymath}
\mathrm{lhs}_1 = (-1)^{\theta_1} \varphi_2(a_{ij} \otimes x_1x_2\cdots x_p a_{ij}^{m-1}, \Delta(y_1y_2\cdots y_q b_{ji}^n)).
\end{displaymath}
For $1\leq j \leq q$, there exists a unique $z_j \in H^-$ such that $\varphi(1,z_j) = 1$ and each of the first tensor factor of $\Delta(y_j) - z_j \otimes y_j$ is  an element of $Y$ strictly greater than $y_j$ multiplied by an element of $H^-$.
 By Lemma \ref{lem: technical lemma} (B), the $\Delta(y_j) - z_j \otimes y_j$ do not contribute to $\mathrm{lhs}_1$. Similarly, for the $n$ copies of $\Delta(b_{ji})$, only one of them contributes $b_{ji} \otimes 1$ to $\mathrm{lhs}_1$, and the rest of them $z \otimes b_{ji}$ with $z = t_{jj}t_{ii}^{-1}$. 
\begin{multline*}
\mathrm{lhs}_1 = (-1)^{\theta_1} \varphi_2(a_{ij} \otimes x_1x_2\cdots x_p a_{ij}^{m-1}, \\
 \prod_{k=1}^q(z_k\otimes y_k) \sum_{l=1}^n (z\otimes b_{ji})^{l-1}(b_{ji}\otimes 1)(z \otimes b_{ji})^{n-l} ).  
\end{multline*}
Note that $\varphi(a_{ij}, z_1z_2\cdots z_q b_{ji}) = \varphi(a_{ij}, b_{ji})$. Also, by Lemma \ref{cor: weight}, 
\begin{displaymath}
b_{ji} z = z b_{ji} \varphi(s_{ii},t_{jj})^{-1}\varphi(s_{jj},t_{ii})^{-1} \varphi(s_{ii},t_{ii})\varphi(s_{jj},t_{jj}) = z b_{ji} \tau_{ij} (-1)^{|b_{ji}|}. 
\end{displaymath}
Thus $(z\otimes b_{ji})^{l-1}(b_{ji}\otimes 1)(z \otimes b_{ji})^{n-l} = (-1)^{(n-1)|b_{ji}|} \tau_{ij}^{n-l} z^{n-1}b_{ji} \otimes b_{ji}^{n-1}$ and
\begin{align*}
\mathrm{lhs}_1 &= (-1)^{\theta_1+\theta_2} (n)_{\tau_{ij}} \varphi_2(a_{ij} \otimes x_1x_2\cdots x_p a_{ij}^{m-1}, b_{ji} \otimes y_1y_2\cdots y_q b_{ji}^{n-1}) \\
&= (n)_{\tau_{ij}}  \varphi_2(x_1x_2\cdots x_p a_{ij}^{m-1} \otimes a_{ij}, y_1y_2\cdots y_q b_{ji}^{n-1} \otimes b_{ji}).
\end{align*}
Here $\theta_2 = |b_{ji}||b_{ji}^{n-1}y_1y_2\cdots y_q|$. In the second identity observe that $\varphi$ respects the parity: $\varphi_2(a\otimes b,c\otimes d) = \varphi_2(b\otimes a,d\otimes c) \times (-1)^{|a||b|+|c||d|}$. The rest is clear from the induction hypothesis. \hfill $\Box$

Let $\Gamma$ be the set of functions $f: X \longrightarrow \BZ_{\geq 0}$ such that $f(x) \leq 1$ if $|x| = \odd$. Such an $f$ induces, by abuse of language, another function $f: Y \longrightarrow \BZ_{\geq 0}$ defined by $f(b_{ji}) := f(a_{ij})$. Set
\begin{equation}  \label{equ: PBW bases}
a_f := \prod_{x \in X}^{\succ} x^{f(x)} \in U^>,\quad b_f := \prod_{y\in Y}^{\succ} y^{f(y)} \in U^<.
\end{equation}
Here $\prod\limits^{\succ}$ means the product with descending order. If $i \leq M < j$, then $\tau_{ij} = -1$ and $\varphi(a_{ij}^m,b_{ji}^m) = 0$ for $m > 1$, which is the reason for $f(a_{ij}) \leq 1$.
\begin{corollary}  \label{cor: orthogonal property}
For $f,g \in \Gamma$ we have $\varphi(a_f,b_g) \neq 0$ if and only if $f = g$. Moreover, the $a_f$ and the $b_f$ form  bases of $U^>$ and $U^<$ respectively. 
\end{corollary}
\proof
The first statement comes from Lemmas \ref{lem: technical lemma}--\ref{cor: technical corollary}; notably the $a_f$ (resp. the $b_f$) are linearly independent. For the second statement, consider $U^>$ for example. A slight modification of the arguments in the proof of \cite[Lemma 2.1]{MRS} by using $R_{23}S_{12}S_{13} = S_{13}S_{12}R_{23}$ shows that $U^>$ is spanned by ordered products of the $a_{ij}$. It remains to prove $s_{ij}^2 = 0$ (and so $a_{ij}^2 = 0$) if $s_{ij}$ is odd; this comes from a comparison of coefficients of $v_i \otimes v_i$ in the equality $R_{23}S_{12}S_{13}(v_j \otimes v_j) = S_{13}S_{12}R_{23}(v_j\otimes v_j) \in U \otimes \BV^{\otimes 2}$.  
\hfill $\Box$
\section{Universal R-matrix}   \label{sec: universal r matrix}
In this section we compute the universal R-matrix of $U_{r,s}$. For this purpose, we first work with a topological version of quantum supergroups and view $r,s$ as formal variables: 
\begin{displaymath}
r = e^{\hbar} \in \BC[[\hbar,\wp]],\quad s = e^{\wp} \in \BC[[\hbar,\wp]].
\end{displaymath} 

\noindent {\bf Step 1.} Extend $U^{\pm},U$ to topological Hopf superalgebras over $\BC[[\hbar,\wp]]$ based on the weight grading : first add commutative primitive elements $(\epsilon_i^*)_{1\leq i \leq M+N}$ of even parity such that $[\epsilon_i^*,x] = \lambda_i x$ for $x \in U^{\pm},U$ of weight $\lambda = \sum_i \lambda_i \epsilon_i \in \BP$; then identify (for the indexes $1\leq i,j \leq M+N$)
\begin{equation*}
s_{ii} = e^{(\hbar + \wp)\sum\limits_{j<i} \epsilon_j^*} \times \begin{cases}
e^{\hbar \epsilon_i^* } & (i \leq M), \\
e^{\wp \epsilon_i^*} & (i > M),
\end{cases} \quad t_{ii} = e^{(\hbar + \wp)\sum\limits_{j<i} \epsilon_j^*} \times \begin{cases}
e^{\wp \epsilon_i^* } & (i \leq M), \\
e^{\hbar \epsilon_i^*} & (i > M).
\end{cases} 
\end{equation*}
Denote by $U^{\pm}_{\hbar,\wp},U_{\hbar,\wp}$ the resulting topological Hopf superalgebras. Set
$$U_{\hbar,\wp}^{\pm} \ni H_i := (\hbar + \wp)\sum_{j<i}\epsilon_j^* +\epsilon_i^* \times \begin{cases}
\wp & (i \leq M), \\
\hbar & (i > M).
\end{cases} $$
 Extend $\varphi$  to a Hopf pairing $\overline{\varphi}: U_{\hbar,\wp}^+ \times U_{\hbar,\wp}^- \longrightarrow \BC((\hbar,\wp))$ by $\overline{\varphi}(\epsilon_i^*,H_j) = \delta_{ij}$. Observe that $\overline{\varphi}(s_{ii},t_{jj}) = \varphi(s_{ii},t_{jj})$, which shows in turn that $\overline{\varphi}$ exists uniquely. The multiplication map induces a surjective morphism of topological Hopf superalgebras from the quantum double $U_{\hbar,\wp}^+ \otimes U_{\hbar,\wp}^-$ to $U_{\hbar,\wp}$ with kernel generated by the $\epsilon_i^* \otimes 1 - 1 \otimes \epsilon_i^*$. 

\noindent {\bf Step 2.} Let $U^0$ be the topological subalgebra of $U_{\hbar,\wp}^{\pm}$ generated by the $\epsilon_i^*$. Then $ U_{\hbar,\wp}^{+} = U^0 U^>$ and $U_{\hbar,\wp}^- = U^0 U^<$.  Corollary \ref{cor: orthogonal property} still holds true. We obtain orthonormal bases of $\overline{\varphi}$ and the universal R-matrix of $U_{\hbar,\wp}$:
\begin{align}
\mathcal{R} &:= \mathcal{R}^0  \mathcal{R}^+,\quad \mathcal{R}^+ =  \sum_{f \in \Gamma} (-1)^{|a_f|} \frac{a_f \otimes b_f}{\varphi(a_f,b_f)} = \prod_{i<j}^{\succ} \mathcal{R}_{ij},  \label{equ: universal R matrix factorization}   \\
\mathcal{R}^0 &= \prod_i e^{\epsilon_i^* \otimes H_i} = \prod_{i\leq M} s^{\epsilon_i^* \otimes \epsilon_i^*} \times \prod_{j>M} r^{\epsilon_j^* \otimes  \epsilon_j^*} \times \prod_{l < k} (rs)^{\epsilon_k^* \otimes \epsilon_l^*},  \label{equ: R matrix Cartan part}   \\
\mathcal{R}_{ij} &=  \begin{cases}
\sum\limits_{n= 0}^{\infty} \frac{a_{ij}^n \otimes b_{ji}^n}{(n)_{rs^{-1}}^! (s^{-1}-r^{-1})^n} & \textrm{if}\ (i<j\leq M), \\
\sum\limits_{n= 0}^{\infty} \frac{a_{ij}^n \otimes b_{ji}^n}{(n)_{sr^{-1}}^! (r^{-1}-s^{-1})^n} & \textrm{if}\ (M<i<j), \\
1 - \frac{a_{ij} \otimes b_{ji}}{s^{-1}-r^{-1}} & \textrm{if}\ (i \leq M < j).
\end{cases}   \label{equ: partial R matrix}
\end{align} 
The formula of $\mathcal{R}$ is similar to that for $U_q(\mathfrak{gl}(M,N))$ in \cite[\S 10.6]{Y1} when $r = q = s^{-1}$. We shall evaluate $\mathcal{R}$ in certain representations (defined over $\BC$).  

\noindent {\bf Step 3.} Assume that $r,s \in \BC^{\times}$ and $\frac{r}{s}$ is not a root of unity. We work with $U_{r,s} = U$ instead of $U_{\hbar,\wp}$.  Let $V$ be a $U$-module (over $\BC$) and $\lambda = \sum_i \lambda_i \epsilon_i \in \BP$. Define $V_{\lambda}$ to be the subspace of $V$ formed of vectors $v$ such that:
 \begin{equation}  \label{equ: weight vector}
s_{ii} v = \varphi(s_{ii},t_{ii})^{\lambda_i} (rs)^{\sum_{j<i}\lambda_j} v,\quad t_{ii}v = \varphi(s_{ii},t_{ii})^{-\lambda_i} (rs)^{\sum_{j\leq i}\lambda_j} v
\end{equation}
for all $1\leq i \leq M+N$.
If $V_{\lambda} \neq 0$, then it is called a weight space of weight $\lambda$. By Lemma \ref{cor: weight}, if $x \in U$ is of weight $\mu$, then $x V_{\lambda} \subseteq V_{\lambda+\mu}$.

 Define $\BQ$ (resp. $\BQ^+$) to be the $\BZ$-span (resp. the $\BZ_{\geq 0}$-span) of the $\epsilon_i - \epsilon_j$ for $i < j$. As in \cite{BW1}, $V$ is said to be in category $\mathcal{O}$ if: (i)  it is spanned by weight spaces; (ii) all the weight spaces are  finite-dimensional; (iii) the set of weights is contained in  
$\cup_{\lambda \in F} (\lambda - \BQ^+)$ for some finite subset $F \subset \BP$.

Let $V, W$ be in category $\mathcal{O}$. Then $\mathcal{R}^0_{V,W} \in \Hom(V \otimes W)$ is well-defined:\footnote{In the non-graded case $\mathcal{R}^0_{V,W}$ is exactly the operator $s \times \widetilde{f}_{V,W}$ in \cite[\S 4]{BW1}.}
$$  v \otimes w \mapsto v \otimes w \times s^{\sum_{i\leq M} \lambda_i \mu_i} r^{\sum_{j>M} \lambda_j \mu_j} (rs)^{\sum_{k>l} \lambda_k\mu_l}  $$
for $v \in V_{\lambda}$ and $w \in W_{\mu}$ where $\lambda = \sum_i \lambda_i \epsilon_i$ and $\mu = \sum_i \mu_i \epsilon_i$. Next, for $f \in \Gamma$, the weight of $a_f v \in V$ is $\lambda + \sum_{i<j} f(a_{ij}) (\epsilon_i-\epsilon_j)$. By condition (iii), $a_f v = 0$ for all but finitely many $f$. So $\mathcal{R}^+_{V,W}\in \Hom(V \otimes W)$ is indeed a finite sum. Let $\mathcal{R}_{V,W} := \mathcal{R}^0_{V,W} \mathcal{R}^+_{V,W}$.  
From the quantum double construction of $U$ we obtain: Category $\mathcal{O}$ together with the $\mathcal{R}_{V,W}$ is braided. 

Consider the vector representation \eqref{equ: vector repre}. From the proof of Lemma \ref{cor: weight} we see that $v_i$ is of weight $\epsilon_i$ and $\BV$ is in category $\mathcal{O}$. Similar to \cite[\S 10.7]{Y1}:
$$\mathcal{R}_{\BV,\BV} = c_{\BV,\BV} R_{s^{-1},r^{-1}}^{-1} c_{\BV,\BV}. $$

Following \cite{DF,JL}, define Drinfeld--Jimbo generators for $1\leq i < M+N$:
\begin{equation*}   
e_i : = s_{ii}^{-1}s_{i,i+1},\quad f_i :=  t_{i+1,i} t_{ii}^{-1},\quad k_i := s_{ii}^{-1}s_{i+1,i+1}, \quad l_i := t_{i+1,i+1}t_{ii}^{-1}.
\end{equation*}
The following relations are proved in the same way as \cite{JL}:
\begin{align*}
&\Delta(e_i) = 1 \otimes e_i + e_i \otimes k_i,\quad \Delta(f_j) = l_j \otimes f_j + f_j \otimes 1, \\
&e_i^2e_{i+1} - (r+s)e_ie_{i+1}e_i + rs e_{i+1}e_i^2 = 0 \quad \textrm{if}\ (1\leq i < M+N-1, i \neq M),   \\
&e_{i-1}e_i^2 - (r+s)e_ie_{i-1}e_i + rs e_i^2e_{i-1} = 0 \quad \textrm{if}\ (1<i<M+N, i \neq M),   \\
&rsf_i^2f_{i+1} - (r+s)f_if_{i+1}f_i + f_{i+1}f_i^2 = 0 \quad \textrm{if}\ (1\leq i < M+N-1, i \neq M),  \\
&rsf_{i-1}f_i^2 - (r+s)f_if_{i-1}f_i + f_i^2f_{i-1} = 0 \quad \textrm{if}\ (1<i<M+N, i \neq M),   \\
&e_i e_j = e_je_i, \quad f_i f_j = f_jf_i, \quad e_M^2 = f_M^2 = 0 \quad \textrm{if}\ (|i-j|>1),  \\
&[e_i, f_j] = \delta_{ij} (-1)^{|i|} (s^{-1}-r^{-1})(k_i - l_i), \\
& e_{M-1}e_Me_{M+1}e_M + rs e_{M+1}e_Me_{M-1}e_M + e_Me_{M-1}e_Me_{M+1}  \\
& \quad \quad +  rse_{M}e_{M+1}e_Me_{M-1} - (r+s) e_Me_{M-1}e_{M+2}e_M = 0 \quad \textrm{if}\ M, N > 1, \\
& rs f_{M-1}f_Mf_{M+1}f_M + f_{M+1}f_Mf_{M-1}f_M + rs f_Mf_{M-1}f_Mf_{M+1}  \\
&\quad \quad + f_{M}f_{M+1}f_Mf_{M-1} - (r+s) f_Mf_{M-1}f_{M+2}f_M = 0 \quad \textrm{if}\ M, N > 1. 
\end{align*}

Let $R' := c_{\BV,\BV} R^{-1}c_{\BV,\BV}$. Then $rsR'$ in the non-graded case is the R-matrix \cite[Definition 3.1]{JL} defining the two-parameter quantum group. The generators $l_{ij}^+$ and $l_{ji}^-$ therein correspond to our $s_{ij}$ and $t_{ji}$.

\end{document}